\documentclass[12pt]{amsart}
\usepackage{amsmath,amssymb}   

\newtheorem{theorem}{Theorem}[section]                    
\newtheorem{proposition}[theorem]{Proposition}            
\newtheorem{corollary}[theorem]{Corollary}                
\newtheorem{lemma}[theorem]{Lemma} 

\begin{document}

\title{Free Product Actions and Their Applications.}

\author{Yoshimichi UEDA}
\address{Fuculty of Mathematics, Kyushu University\\ Fukuoka, 810-8560, Japan}
\email{ueda@math.kyushu-u.ac.jp}
\maketitle
\baselineskip=15pt

\begin{abstract}
A quick review of a series of 
our recent works \cite{ueda:minimalaction,ueda:index4,dima-ueda,hiai-ueda} 
on free product actions, partly in collaborations with Dimitri Shlyakhtenko 
and with Fumio Hiai, is given. We also work out type III theoretic subfactor 
analysis on the subfactors constructed in \cite{ueda:index4} (see also 
\cite{dima-ueda}). This part can be read as a supplementary appendix to those works. 
\end{abstract}

\section{Introduction.}

{\it Free product actions} mean group or group-like symmetries on operator 
algebras ($C^*$-algebras or von Neumann algebras), which are naturally arisen 
from the free product construction. In the early 80s, D.~Voiculescu started his 
original study on the II$_1$-factors associated with free groups ${\mathbb F}_n$ 
(so-called free group factors denoted by $L({\mathbb F}_n)$) by the use of idea 
of probability. This line of study is now called {\it free probability theory} (see 
\cite{voiculescu:icm,voiculescu:lecture}), which has been widely accepted as 
one of the major (non-commutative) probability theories, and the free product 
construction plays particular r\^{o}le in the theory, in particular, naturally produces 
many explicit examples of non-commutative probability spaces where the theory 
works quite effectively.    

\medskip
In \cite{ueda:minimalaction}, we began to use free product actions to deal with 
purely operator algebraic questions by the aid of free probability theory, and then, 
with D.~Shlyakhtenko \cite{dima-ueda} (also see \cite{ueda:index4}) we showed 
that any positive real number greater than 4  can be realized as the Jones index 
of an irreducible inclusion of II$_1$-factors both isomorphic to $L({\mathbb F}_{\infty})$. 
We also use, with F.~Hiai \cite{hiai-ueda}, free product actions to construct variety 
of examples of (aperiodic) automorphisms on interpolated free group factors 
$L({\mathbb F}_r)$, $1 < r \leq \infty$, (independently introduced by 
K.~Dykema \cite{dykema}, F.~R{\u{a}}dulescu \cite{radulescu1}).  

\medskip
One of the purposes of this paper is to give a quick review of the above-mentioned 
progress and the other is to give supplementary results to the works 
\cite{ueda:index4,dima-ueda}, part of which was already commented in Remark 3.5 
of \cite{dima-ueda}. The organization is as follows. The \S\ref{SUq},\S\ref{freeproductactions},
\S\ref{wassermannconstruction}  are all preliminaries; namely we recall 
the compact quantum group ${\rm SU}_q(2)$ in \S\ref{SUq}; the notion of 
free product actions in \S\ref{freeproductactions}; 
and Wassermann's construction for minimal actions of compact quantum groups in 
\S\ref{wassermannconstruction}. The \S\ref{review} is quick review of 
\cite{dima-ueda,hiai-ueda}, while in \S\ref{typeIIIsubfactoranalysis} 
we give type III theoretic subfactor analysis on the subfactors constructed in 
\cite{ueda:index4} (also see \cite{dima-ueda}). Finally, in \S\ref{comments} 
we collect some questions and comments related to the subject matter treated 
in this paper.         

The author would like to express his sincere gratitude to Fumio Hiai and Dima 
Shlyakhtenko for their friendship, many discussions and their collaborations with 
the author. The author also thanks Nobuaki Obata, Akito Hora and Taku Matsui for 
giving him this opportunity.  

\section{The quantum group ${\rm SU}_q(2)$.}
\label{SUq} 

The von Neumann algebraic approach to compact quantum groups is needed for our aim, and 
we would like to employ the following as the (rather working) definition of compact quantum 
groups, which is sufficient for the later use. Let $A$ be a von Neumann algebra and $\delta : 
A \rightarrow A\otimes A$ be a (normal) $*$-isomorphism (called a co-multiplication) satisfying 
the co-associativity: $(\delta\otimes\text{Id}_A)\circ\delta = (\text{Id}_A\otimes\delta)\circ\delta$. 
The pair $(A, \delta)$ is usually called a Hopf - von Neumann algebra. A Haar state $h$ 
for $(A,\delta)$ is a faithful normal state on $A$ satisfying the left and right-invariance 
property: 
$(h\otimes\text{Id}_A)\circ\delta(a) = h(a)1_A = (\text{Id}_A\otimes h)\circ\delta(a)$ 
for $a \in A$. We will here call a Hopf - von Neumann algebra equipped with a Haar state  
a (von Neumann algebraic) compact quantum group. 

\medskip
All the usual notions in the group representation theory can be used even in this setting. 
We briefly explain those for the later use. (See \cite{woro2} for details, also see \cite{b-s}.) 
A corepresentation $\pi$ on a finite dimensional Hilbert space $V_{\pi}$ means a linear map 
$\pi : V_{\pi} \rightarrow V_{\pi}\otimes A$ with the property: 
$$
(\pi\otimes\text{Id}_A)\circ\pi = (\text{Id}_{V_{\pi}}\otimes\delta)\circ\pi.
$$
For a basis $\{\xi_j\}$ of $V_{\pi}$ one finds $u(\pi)_{ij}$ in $A$ with    
$$
\pi(\xi_j) = \sum_i \xi_i \otimes u(\pi)_{ij}.    
$$
Via the natural identification $M_{\text{dim}V_{\pi}}(A) = B(V_{\pi})\otimes A$ (with respect 
to the basis $\{\xi_j\}$), the operator matrix $u(\pi) := [u(\pi)_{ij}]$ can be regarded as 
an element in $B(V_{\pi})\otimes A$, and the $u(\pi)_{ij}$'s are called the matrix elements 
with respect to the basis $\{\xi_j\}$. A corepresentation $\pi$ is said to be unitary when the 
operator $u(\pi)$ defines a unitary element. An intertwiner of two corepresentations $\pi_1$, 
$\pi_2$ means a linear map $T \in B(V_{\pi_1},V_{\pi_2})$ satisfying 
$$
(T\otimes1)u(\pi_1) = u(\pi_2)(T\otimes1),
$$ 
and the set of those intertwiners is denoted by $\text{End}(\pi_1,\pi_2)$. When 
$\text{End}(\pi_1,\pi_2) = \{0\}$, we say that $\pi_1$ and $\pi_2$ are nonequivalent. 
A corepresentation $\pi$ is irreducible when the self-intertwiners $\text{End}(\pi) = 
\text{End}(\pi,\pi)$ is 1-dimensional. 

\medskip
Since we consider only examples of compact-type, any finite dimensional corepresentation 
is reducible, that is, it is equivalent to a direct sum of irreducible ones, and also is equivalent 
to, or can be (by replacing the original inner product by a suitable equivalent one), a unitary one.  
Thus, in any case, for a given finite dimensional corepresentation $\pi$, we may and do 
assume that any given $\pi$ is a direct sum of finitely many corepresentations:  
$$
\pi = m_1 \pi_1 \oplus \cdots \oplus m_{n} \pi_n
$$
with mutually nonequivalent irreducible unitary corepresentations $\pi_1, \pi_2,\dots,\pi_n$ 
and their multiplicities $m_1,\dots,m_n$. In this case, one easily describe the 
self-intertwiners $\text{End}(\pi)$ as 
$$
\left(M_{m_1}({\bold C})\otimes{\bold C}1_{V_{\pi_1}}\right)\oplus\cdots\oplus
\left(M_{m_n}({\bold C})\otimes{\bold C}1_{V_{\pi_n}}\right). 
$$

\medskip
The quantum group ${\rm SU}_q(2)$ ($q \in (0,1)$) is a typical example of compact 
quantum group, which is given in the following way. Let us assume the indeteminates 
$x,u,v,y$ satisfy the relations: 
\begin{gather}
ux = q xu,\ vx = q xv,\ yu = q uy,\ yv = q vy, \notag\\ 
uv = vu,\ x^* = y,\ u^* = -q^{-1}v,  \notag
\end{gather}
and consider the universal $C^*$-algebra ${\mathcal A}$ generated by $x,u,v,y$. One 
can then construct a co-multiplication 
$\delta : {\mathcal A} \rightarrow {\mathcal A}\otimes_{\text{min}}{\mathcal A}$ 
by the matrix-product rule (with replacing the usual scalar multiplication by the operator 
tensor product): 
$$
\delta(x) = x \otimes x + u \otimes v, \ 
\delta(u) = x \otimes u + u \otimes y, \
\delta(v) = v \otimes x + y \otimes v, \ 
\delta(y) = v \otimes u + y \otimes y. 
$$
It is known (see \cite{woro1,woro2}) that there is a unique (faithful) Haar state 
(i.e., a state with the left and right-invariance property as above) for the pair 
$({\mathcal A}, \delta)$, and let $A$ be the von Neumann algebra obtained from 
${\mathcal A}$ via the GNS representation of the Haar state. It is straightforward to 
see that the co-multiplication has the unique extension to $A$, and we still denote it 
by $\delta : A \rightarrow A\otimes A$. By the construction, the GNS-vector gives, 
as a vector state, a Haar state $h$ for the pair $(A,\delta)$. This compact quantum 
group $(A,\delta,h)$ is called the quantum group ${\rm SU}_q(2)$, and sometimes 
denote $A = L^{\infty}({\rm SU}_q(2))$ (and ${\mathcal A} = C({\rm SU}_q(2))$).

\section{Free product actions.}
\label{freeproductactions} 

Let $N_1$, $N_2$ be ($\sigma$-finite) von Neumann algebras with faithful normal 
states $\varphi_1$, $\varphi_2$, respectively. We can then construct, from the given 
pairs $(N_1,\varphi_1)$, $(N_2,\varphi_2)$, a new pair $(M,\varphi)$ of a von 
Neumann algebra and a faithful state in such a way that 
\begin{enumerate}\renewcommand{\labelenumi}{(\roman{enumi})}
\item there are two embeddings $N_1, N_2 \hookrightarrow M$ whose ranges 
           generate $M$, i.e., $M = N_1 \vee N_2$; 
\item the restrictions of $\varphi$ to $N_1$, $N_2$ (via the embeddings above) 
           coincide with $\varphi_1$, $\varphi_2$, respectively; 
\item (the {\it free independence}) the expectation value of any alternating word 
           in $N_1^{\circ} := \text{Ker}\varphi_1$, $N_2^{\circ} := \text{ker}\varphi_2$ 
           is $0$. 
\end{enumerate}
The pair $(M, \varphi)$ is called the (reduced) free product of $(N_1,\varphi_1)$  
and $(N_2,\varphi_2)$, and written as  
$$
(M, \varphi) = (N_1, \varphi_1) * (N_2, \varphi_2). 
$$  

The algebra $M$ has many automorphisms naturally induced from automorphisms 
on $N_1$, $N_2$ preserving the states $\varphi_1$, $\varphi_2$, respectively: 
$$
(\alpha_1,\alpha_2) \in \text{Aut}_{\varphi_1}(N_1) \times \text{Aut}_{\varphi_2}(N_2) 
\mapsto \alpha := \alpha_1 * \alpha_2 \in \text{Aut}_{\varphi}(M)
$$
determined (thanks to the characterization of $(M,\varphi)$ explained just above) 
by the relations $\alpha|_{N_1} = \alpha_1$, $\alpha|_{N_2} = \alpha_2$.  We would 
like to refer to group or group-like symmetries of this type as {\it free product actions}, 
and thus the most typical example is the modular action of ${\bf R}$ associated with the 
free state $\varphi$: 
$$
t \in {\bf R} \mapsto \sigma_t^{\varphi}  
= \sigma_t^{\varphi_1} * \sigma_t^{\varphi_2} \in \text{Aut}_{\varphi}(M).
$$ 
This is, of course, a special example, but the reader should pay their attention to it 
because many observations and/or motivations of 
\cite{ueda:minimalaction,ueda:index4,dima-ueda,hiai-ueda} 
came in part from study on modular actions of this type so that many discussions there and/or here 
are closely related to 
the others \cite{ueda:amalgam,ueda:freeproduct,ueda:amalgam-example,ueda:chi-group}. 
For example, the attempt of showing the minimality of free product actions is essentially same game as 
that of finding  the III$_{\lambda}$-types of (type III) free product von Neumann algebras. This is 
quite similar to the infinite tensor product situation (see p.212 of \cite{wassermann} 
for product-type actions of compact Lie groups and p.65--66 of Takesaki's monograph \cite{takesaki} 
for the III$_{\lambda}$-classification of Powers factors).  
 
In the quantum group case, an action does no longer mean an automorphic one, that is,   
an action of a Hopf - von Neumann algebra $(A, \delta)$ on a von Neumann algebra $N$ is 
formulated as a normal $*$-isomorphism $\Psi : N \rightarrow N\otimes A$ satisfying 
$(\Psi\otimes\text{Id}_A)\circ\Psi = (\text{Id}_M\otimes\delta)\circ\Psi$. 
Assume that two actions 
$\Gamma_1 : N_1 \rightarrow N_1\otimes A$, 
$\Gamma_2 : N_2 \rightarrow N_2\otimes A$ 
of a Hopf - von Neumann algebra $(A,\delta)$ on von Neumann algebras $N_1$, $N_2$ 
equipped with faithful normal states $\varphi_1$, $\varphi_2$, respectively, satisfy 
the following invariance condition:   
$$
(\varphi_1\otimes\text{Id}_A)\circ\Gamma_1(x) = \varphi_1(x), \quad 
(\varphi_2\otimes\text{Id}_A)\circ\Gamma_2(y) = \varphi_1(y)
$$
for $x \in N_1$, $y \in N_2$. Then one can, even in this case, construct the free product action 
$$
\Gamma := \Gamma_1 * \Gamma_2 : M \rightarrow M\otimes A
$$
of $(A,\delta)$ on the free product 
$$
(M,\varphi) = (N_1,\varphi_1) * (N_2,\varphi_2)
$$ 
(by the essentially same way as in the group case) in such a way that 
$$
\Gamma|_{N_1} = \Gamma_1, \quad \Gamma|_{N_2} = \Gamma_2. 
$$
See \cite{ueda:minimalaction,ueda:index4} for details.

\section{Wassermann's construction for minimal actions of compact quantum groups.} 
\label{wassermannconstruction}

Wassermann's construction (see \cite{wassermann}) is a way of constructing 
subfactors of finite-index from (minimal) actions of compact groups together with their 
(finite dimensional unitary) representations. In this section, we will explain how it can 
be generalized to the quantum group case, which was done in \cite{ueda:index4}. 
(See also \cite{banica2}, where the only Kac case is treated.) For notational simplicity, 
we will deal with only the ${\rm SU}_q(2)$ (i.e., $(A := L^{\infty}({\rm SU}_q(2)),\delta,h)$, 
just explained in \S\ref{SUq}), but this is not essential; in fact we will not use any special 
characteristic of $\text{\rm SU}_q(2)$.  

\medskip
Let $\Gamma : M \rightarrow M\otimes A$ be an action of ${\rm SU}_q(2)$ on 
a factor $M$. Let us assume that 
\begin{enumerate}\renewcommand{\labelenumi}{(\alph{enumi})}
\item the fixed-point algebra $M^{\Gamma} := \{ x \in M : \Gamma(m) = m\otimes 1\}$ 
has the trivial relative commutant  in $M$; 
\item for any (irreducible) unitary corepresentation $(\pi, V_{\pi})$ of $SU_q(2)$ one can 
find a Hilbert space ${\mathcal H}$ in $M$ in such a way that $(\Gamma|_{\mathcal H}, 
{\mathcal H})$ is unitarily equivalent to $(\pi, V_{\pi})$. See \cite{roberts} for 
the terminology of {\it Hilbert spaces in von Neumann algebras}.
\end{enumerate}
The actions given in \cite{ueda:minimalaction,ueda:index4,dima-ueda} 
clearly satisfy these conditions by the constructions themselves, but in general the latter (b) 
should be replaced by the faithfulness in the sense of Izumi-Longo-Popa \cite{izumi-longo-popa} 
with the aid of some (standard) techniques along the exactly same line as in \cite{akth,roberts}. 
The reader is encouraged to consult with Nakagami-Takesaki's monograph \cite{nakagami-takesaki} 
although only group actions are treated there. In fact, it is quite a straightforward work to 
generalize the necessary parts given there to the quantum group case; see 
\cite{boca,izumi-longo-popa,yamanouchi} for example. 

\medskip 
Let $\pi$ be a finite dimensional corepresentation on $V_{\pi}$ of ${\rm SU}_q(2)$. 
We then consider the following inclusion: 
$$
{\mathcal M} := (B(V_{\pi})\otimes M)^{\text{Ad}\pi\otimes\Gamma} \supseteq 
{\mathcal N} := {\bold C}1_{V_{\pi}}\otimes M^{\Gamma}, 
$$
where $\text{Ad}\pi\otimes\Gamma : B(V_{\pi})\otimes M \rightarrow 
B(V_{\pi})\otimes M\otimes A$ is the action defined by 
$$
\text{Ad}\pi\otimes\Gamma := 
\text{Ad}u(\pi)_{13}\circ(\text{Id}_{B(V_{\pi})}\otimes\Gamma).
$$
See \S\ref{SUq} for the definition of $u(\pi)$. Thanks to the assumption (b) we can construct  
an explicit isomorphism between ${\mathcal M}$ and ${\mathcal N}$, and hence, 
in particular, ${\mathcal M}$ is a factor (see Remark 10 of \cite{ueda:index4} or 
Theorem 3.1 of \cite{dima-ueda}). When replacing the (b) by the faithfulness, 
one should do the argument after taking tensor product with $B({\mathcal H})$, and 
then remove it so that the consequence is slightly changed, i.e., ${\mathcal M}$ and 
${\mathcal N}$ are stably isomorphic.

\medskip
The first non-trivial problem is whether there is a (finite-index) conditional expectation 
for ${\mathcal M} \supseteq {\mathcal N}$. This is somewhat non-trivial because of 
the definition of the action $\text{Ad}\pi\otimes\Gamma$ and the non-commutativity 
of the algebra $A$. However, we can indeed overcome this by showing the following:

\begin{lemma} {\rm (See p.44--45 of \cite{ueda:index4}.)} 
Let $\{f_z\}_{z \in {\bold C}}$ be the Woronowicz characters of ${\rm SU}_q(2)$ 
{\rm (}see {\rm \cite{woro2})}, and consider the element 
$F_{\pi} \in B(V_{\pi})$ constructed as $F_{\pi}^z := (\text{\rm Id}_{B(V_{\pi})}
\otimes f_z)(u(\pi))$, $z \in {\bold C}$. Let us consider the quantum dimension and 
the positive / negative $q$-traces 
$$
\text{\rm dim}_q(\pi) := \text{\rm Tr}(F_{\pi}) = \text{\rm Tr}(F_{\pi}^{-1}), 
$$$$
\tau_q^{(\pi,+)}(\ \cdot\ ) := 
\frac{1}{\text{\rm dim}_q(\pi)}\text{\rm Tr}(F_{\pi}\ \cdot\ ), \quad 
\tau_q^{(\pi,-)}(\ \cdot\ ) := 
\frac{1}{\text{\rm dim}_q(\pi)}\text{\rm Tr}(F_{\pi}^{-1}\ \cdot\ ), 
$$
where $\text{\rm Tr}$ means the non-normalized usual trace on $B(V_{\pi})$. Then we 
have
$$
(\tau_q^{(\pi,-)}\otimes\text{\rm Id}_M)\circ E_{\text{Ad}\pi\otimes\Gamma} = 
E_{\Gamma}\circ(\tau_q^{(\pi,+)}\otimes\text{\rm Id}_M), 
$$
where $E_{\Gamma}$, $E_{\text{Ad}\pi\otimes\Gamma}$ denote the usual conditional 
expectations onto the fixed-point algebras constructed by using the Haar state $h$. 
\end{lemma}

\medskip
Therefore, the restriction of the Fubini-map $\tau_q^{(\pi,-)}\otimes\text{Id}_M$ to ${\mathcal M}$ 
gives rise to a (faithful normal) conditional expectation $E : {\mathcal M} \rightarrow {\mathcal N}$. 
Unlike in the usual group case, we cannot provide the invariance principle of the exactly same form 
as in \cite{wassermann} because of the same reason as in the trouble treated above, and 
hence we start with constructing a candidate of the triple of basic extension, dual conditional 
expectation and Jones projection, and then show that they satisfy the so-called Pimsner-Popa
characterization. 

\medskip
Consider the bigger algebra 
$$
{\mathcal M}_1 := \left(B(V_{\hat{\pi}}\otimes V_{\pi})\otimes M\right)^
{\text{Ad}(\hat{\pi}\otimes\pi)\otimes\Gamma}
$$
into which the inclusion ${\mathcal M } \supseteq {\mathcal N}$ is embedded with 
$x \mapsto 1_{V_{\hat{\pi}}}\otimes x$, where $\hat{\pi}$ means the {\it canonical dual} 
corepresentation on $V_{\hat{\pi}} := \overline{V}_{\pi}$ of $\pi$ (see \cite{banica1}). 
By the same reason as above, it can be checked that ${\mathcal M}_1$ is a factor, and also then 
that the restriction of the Fubini-map 
$\tau_q^{(\hat{\pi},-)}\otimes\text{Id}_{B(V_{\pi})}\otimes\text{Id}_M$ to ${\mathcal M}_1$ 
gives rise to a conditional expectation $E_1 : {\mathcal M}_1 \rightarrow {\mathcal M}$ 
thanks to the equation 
$\text{Ad}(\hat{\pi}\otimes\pi)\otimes\Gamma = 
\text{Ad}\hat{\pi}\otimes(\text{Ad}\pi\otimes\Gamma)$. 
Looking at the basic extension of the inclusion
$$
B(V_{\pi}) \supseteq {\bold C}1, \quad 
\tau_q^{(\pi,-)}(\ \cdot\ ) : B(V_{\pi}) \rightarrow {\bold C}1,
$$
one can easily show that 
$$
{\mathcal M}_1 \overset{E_1}{\supseteq} {\mathcal M} \overset{E}{\supseteq} {\mathcal N}
$$
satisfies one variation (see Theorem 8 of \cite{hamachi-kosaki3}) of 
the Pimsner-Popa characterization of basic extensions so that this tower is the Jones tower 
up to the $1$st step. As a consequence one sees that the index $\text{Ind}E$ is given by 
the square of the quantum dimension $(\text{dim}_q(\pi))^2$. 

\medskip
Summing up the above discussions we in particular have the following: 

\begin{proposition} {\rm (\cite{ueda:index4})} We have \newline
\ \ {\rm (1)} ${\mathcal M} \cong {\mathcal N}$, and in particular, both being factors; \newline
\ \ {\rm (2)} $\text{\rm Ind}E = (\text{\rm dim}_q\pi)^2 
           := \text{Tr}(F_{\pi})^2 = \text{Tr}(F_{\pi}^{-1})^2$. 
 
If one uses the faithfulness together with the usual technique of taking tensor product 
with $B({\mathcal H})$ instead of the assumption (2), the first assertion should be 
changed into \newline
\ \ {\rm (1')} ${\mathcal M}$ is stably isomorphic to ${\mathcal N}$. \newline 
{\rm (}Of course, this change is not needed when the fixed-point algebra is properly infinite.{\rm )}  
\end{proposition}

\medskip
Clearly, the above procedure of finding the basic extension can be iterated, and the $n$-th basic 
extension ${\mathcal M_n}$ has the simple description: 
$$
{\mathcal M}_n = \left(B(V_{\rho_n})\otimes M\right)^
{\text{Ad}\rho_n\otimes\Gamma}
$$
with the embedding $x \mapsto 1\otimes x$, where $\rho_n$ is determined by 
the recursive relations: 
$$
\begin{aligned}
\rho_{2m + 1} &= \hat{\pi}\otimes\rho_{2m}, \\
\rho_{2m} &= \pi\otimes\rho_{2m-1}, \\
\rho_0 &= \pi. 
\end{aligned}
$$   
We should point out here that the conditional expectation $E_n : {\mathcal M}_n \rightarrow 
{\mathcal M}_{n-1}$ is of course described by the same way as above based on the $q$-traces.  

\section{Review on the works with D.~Shlyakhtenko and with F.~Hiai.} 
\label{review}

\subsection{Irreducible subfactors of $L({\mathbb F}_{\infty})$: Joint work with 
D.~Shlyakhtenko.} 

Based on \cite{ueda:index4} with the aid of free Araki-Woods factors introduced 
by D.~Shlyakhtenko \cite{dima1,dima2} (also see \cite{dima3}), we showed with 
D.~Shlyakhtenko the following theorem: 

\medskip
\begin{theorem}
\label{dima-ueda-maintheorem}
 {\rm (\cite{dima-ueda})} 
Any real number $>4$ can be realized as the Jones index $[{\mathcal M}:{\mathcal N}]$ 
of an irreducible inclusion ${\mathcal M} \supseteq {\mathcal N}$ in such a way that 
${\mathcal M} \cong {\mathcal N} \cong L({\mathbb F}_{\infty})$. 
\end{theorem}

\medskip
This together with his pioneering work of F.~R\u{a}dulescu \cite{radulescu1} in the subject 
matter implies 

\begin{corollary} The possible Jones indices of irreducible inclusions consisting of 
$L({\mathbb F}_{\infty})$ is the so-called Jones spectrum 
$$
\left\{4\cos^2\frac{\pi}{n+1} : n \geq 2 \right\} \cup [4,\infty]. 
$$
\end{corollary}

\medskip
In the proof, we first construct the free product action $\Gamma : M \rightarrow M\otimes A$ 
of ${\rm SU}_q(2)$ by  
$$
(M, \varphi) := (L({\mathbb F}_N), \tau_{{\mathbb F}_N}) * (A, h), \quad 
\Gamma := (\text{Id}_{L({\mathbb F}_N)}\otimes 1_A) * \delta,  
$$
and then consider the Wassermann-type inclusion associated with an irreducible corepresentation. 
This inclusion consists of type III$_{q^2}$ factors so that we show that it is of essentially type II, 
and then get the desired type II$_1$ subfactor via the structure theorem of type III factors. What 
we have to do is to investigate the fixed-point algebra $M^{\Gamma}$, or more precisely its discrete 
core. In the process of identifying with $L({\mathbb F}_{\infty})$, one of the free Araki-Woods factors 
(see \cite{dima1}) naturally arises in the following way: Our type III subfactor $M^{\Gamma}$ 
can be described, via the Takesaki duality (for compact quantum group action), as the amalgamated 
free product 
$$
\left(L({\mathbb F}_N)\otimes\widehat{A}, \tau_{{\mathbb F}_N}\otimes\text{Id}\right) 
*_{\hat{A}} \left(B(L^2(A)), F \right), 
$$
where $\widehat{A} := \ell^{\infty}(\widehat{{\rm SU}_q(2)})$ means the dual Hopf-von Neumann 
algebra of ${\rm SU}_q(2)$. Since $\widehat{A}$ is an infinite direct sum of matrix algebras, 
an argument similar to Theorem 5.1 of \cite{dima2} works for identifying $M^{\Gamma}$ 
(or the amalgamated free product of the above form) with the unique free Araki-Woods factor of 
type III$_{q^2}$. The free Araki-Woods factor of type III$_{\lambda}$ ($\lambda \neq 0,1$) is known 
to have the discrete core isomorphic to $L({\mathbb F}_{\infty})\otimes B({\mathcal H})$, and hence 
we observe an exact relationship between our subfactor and the free group factor 
$L({\mathbb F}_{\infty})$. 

\medskip
Therefore, one of the key-points is ``passing once through infinite (type III) factors to get 
a result on type II$_1$ factors." This should be much emphasized. In fact, S.~Popa and 
D.~Shlyakhtenko \cite{popa-dima} generalized this theorem further to the final form, 
where the use of type II$_{\infty}$ factors is one of the most essential points although 
the counter-part of the arguments explained above is much more complicated in their case. 

\subsection{Automorphisms of $L({\mathbb F}_r)$ with $r > 1$: Joint work with 
F.~Hiai.} 
 
Using free product actions, we exhibited, with F.~Hiai, several automorphisms with 
special characteristics on each (interpolated) free group factor $L({\mathbb F}_r)$ 
with $r > 1$ (\cite{dykema,radulescu1}).  

\medskip
\begin{theorem} {\rm (\cite{hiai-ueda})} 
There are continuously many non-outer conjugate aperiodic automorphisms on 
the interpolated free group factor $L({\mathbb F}_r)$ with any fixed $r > 1$, all of 
whose crossed-products are isomorphic. 
\end{theorem} 

\medskip
This theorem is a supplementary result to a famous one of J.~Phillips (see 
\cite{phillips1,phillips2}) obtained in the 70s. We also worked on this isomorphic 
crossed-product algebra; it has the property $\Gamma$ but not McDuff or not strongly 
stable, and its quotient group $\overline{\rm Int}/{\rm Int}$ is isomorphic to 
the 1st cohomology group of the unique amenable type II$_1$ equivalence relation. 
Those proofs are essentially same as those given in \cite{ueda:chi-group}. 

\medskip
Let $P$, $Q$ be von Neumann algebras equipped with faithful normal states 
$\varphi_P$, $\varphi_Q$, respectively. Let $G$ be a discrete (countable) group, 
and the compact co-commutative Kac algebra ${\mathbb K}_G = (L(G),\delta_G,
\kappa_G,\tau_G)$ can act on the free product algebra 
$$
(M, \varphi) := (P\otimes L(G), \varphi_P\otimes\tau_G) * (Q,\varphi_Q)
$$
by the free product action
$$
\Gamma^G =
(\text{Id}_P\otimes\delta_G) * (\text{Id}_Q\otimes 1_{L(G)}).
$$
This action is related to the following action $\alpha$ of $G$ on 
the free product algebra 
$$
(N, \psi) = (P, \varphi_P) *
\left(\underset{g \in G}{*} (Q, \varphi_Q)_g \right)
$$
defined by 
$$
\alpha_g := \text{Id}_P * \gamma^G_g
$$ 
with the $G$-free shift action $\gamma^G$ on 
$\displaystyle{\underset{g \in G}{*}(Q, \varphi_Q)_g}$. 
In fact, the following duality holds:  

\begin{theorem}\label{hiai-ueda-duality} {\rm (\cite{hiai-ueda})} We have
$$
(M,\Gamma^G) \cong (N\rtimes_{\alpha}G,\widehat{\alpha})
$$ 
with the dual action $\widehat{\alpha}$ 
{\rm ($= \widehat{\text{Id}_P * \gamma^G}$)} 
{\rm (}see {\rm \cite{nakagami-takesaki})},
\end{theorem} 

\medskip
This theorem explains (rather philosophically) that the minimal actions of 
${\rm SU}_q(N)$ constructed in \cite{ueda:minimalaction,ueda:index4} 
can be regarded as the dual actions of ``free Bernoulli shift actions" of the dual 
$\ell^{\infty}(\widehat{{\rm SU}_q(N)})$. We will return to this point of view 
in near future. 

\medskip
This duality can be used for example to investigate the crossed-products by free Bernoulli 
shifts on the free group factor $L({\mathbb F}_{\infty})$, see \cite{hiai-ueda} 
for details. 

\section{Type III theoretic subfactor analysis on Wassermann-type subfactors associated with 
${\rm SU}_q(2)$} 
\label{typeIIIsubfactoranalysis}

Throughout this section, we keep and use freely the notations in 
\S\ref{SUq},
\S\ref{freeproductactions},
\S\ref{wassermannconstruction}. 
Let us assume that $\Gamma : M \rightarrow M\otimes A$ is a minimal action of 
$\text{\rm SU}_q(2)$ equipped with an invariant state $\varphi$ satisfying 
\begin{enumerate}\renewcommand{\labelenumi}{(\roman{enumi})}
\item the fixed-point algebra $M^{\Gamma}$ is of type III$_{q^2}$; 
\item the modular automorphism $\sigma_t^{\varphi}$ has the period 
           $\displaystyle{T_0 := -\frac{2\pi}{\log q^2}}$; 
\item there is a von Neumann subalgebra $Q$ of $M^{\Gamma}\cap M^{\sigma^{\varphi}}$ 
           with $Q'\cap M = {\bold C}1$. 
\end{enumerate}
These conditions are typical characteristic that the actions given in 
\cite{ueda:minimalaction,ueda:index4,dima-ueda} 
(one of which is explained in \S\S5.1 briefly) has. Let us next choose an arbitrary finite 
dimensional unitary corepresentation $\pi : V_{\pi} \rightarrow V_{\pi}\otimes A$, 
and consider the Wassermann-type inclusion consisting of (isomorphic) type III$_{q^2}$ 
factors 
$$
{\mathcal M}({\pi}) := \left(B(V_{\pi})\otimes M\right)^{\text{\rm Ad}\pi\otimes\Gamma} 
\overset{E_{\pi}}{\supseteq} 
{\mathcal N}({\pi}) := {\bold C}1_{V_{\pi}}\otimes M^{\Gamma}
$$ 
Since we do not assume that $\pi$ is irreducible, the inclusion may not be irreducible so that 
it is unclear whether or not the inclusion ${\mathcal M}(\pi) \supseteq {\mathcal N}(\pi)$ 
is of essentially type II. In what follows, we will discuss this aspect completely. Namely, 
following Hamachi-Kosaki's analysis on type III subfactors   
(\cite{hamachi-kosaki1,hamachi-kosaki2,kosaki1,kosaki2,kosaki3}), 
we investigate the factor map structure and the three-step inclusions associated with 
the type III inclusion ${\mathcal M}(\pi) \supseteq {\mathcal N}(\pi)$, and as one consequence 
give a necessary and sufficient condition, in terms of fusion rule of $\pi$ decomposed into 
irreducible corepresentations, when the inclusion in question is of essentially type II, 
i.e., it can be described by an inclusion of type II$_1$ factors (together with a certain 
${\mathbb Z}$-action on its amplification) in principle. We encourage the reader to consult 
with \cite{kosaki2} as a brief survey of Hamachi-Kosaki's analysis.  

\medskip
As explained in \S\ref{SUq}, we may and do assume that the $\pi$ is the direct sum of 
(unitary) irreducible corepresentations:  
$$
\pi = m_1 \pi_1 \oplus \cdots \oplus m_n \pi_n,  
$$
and hence  
$$
\text{\rm End}(\pi) =   
\left(M_{m_1}({\bold C})\otimes{\bold C}1_{V_{\pi_1}}\right)\oplus\cdots\oplus
\left(M_{m_n}({\bold C})\otimes{\bold C}1_{V_{\pi_n}}\right). 
$$
Therefore, one has: 
$$
{\mathcal N}(\pi)' \cap {\mathcal M}(\pi) = 
B(V_{\pi})^{\text{\rm Ad}\pi}\otimes{\bold C}1_M \cong 
B(V_{\pi})^{\text{\rm Ad}\pi} = 
\text{\rm End}(\pi) 
$$
(since $\displaystyle{\left(M^{\Gamma}\right)'\cap M = {\bold C}1_M}$) so that the inclusion 
${\mathcal M}(\pi) \supseteq {\mathcal N}(\pi)$ is irreducible if and only if so is $\pi$. 

\bigskip
Let us start with looking at the inclusion: 
$$
\widetilde{{\mathcal M}(\pi)} := {\mathcal M}(\pi)\rtimes_{\sigma^{\psi\circ E_{\pi}}}{\bold R} 
\quad \overset{\hat{E}_{\pi}}{\supseteq} \quad  
\widetilde{{\mathcal N}(\pi)} := {\mathcal N}(\pi)\rtimes_{\sigma^{\psi}}{\bold R}, 
$$
where $\psi$ is the restriction of $\varphi$ to $M^{\Gamma} = {\mathcal N}(\pi)$. 
Here $\hat{E}_{\pi}$ is defined by  
$$
\hat{E}_{\pi}\left(\int_{-\infty}^{\infty} m(t)\lambda(t) dt\right) := 
\int_{-\infty}^{\infty} E_{\pi}(m(t))\lambda(t) dt.
$$

\medskip
The following simple application of the Fourier transform is essential in our computations: 

\begin{lemma}\label{haagerupstormerlemma}
{\rm (e.g. Proposition 5.6 of \cite{haagerup-stormer})} 
There is a surjective isomorphism 
$$
\Phi :  M^{\Gamma}\rtimes_{\sigma}{\bold R} \rightarrow 
\left(M^{\Gamma}\rtimes_{\dot{\sigma}}{\mathbb T}\right)
\otimes L^{\infty}[0,-\log q^2)
$$
determined by 
$$
\Phi\left(\pi_{\sigma}(m)\right) = \pi_{\dot{\sigma}}(m)\otimes 1, \quad 
\Phi\left(\lambda^{\sigma}(t)\right) = 
\lambda^{\dot{\sigma}}([t])\otimes m_{e^{it(\ \cdot\ )}}. 
$$
Here, we set $\sigma := \sigma^{\psi}$ and $\dot{\sigma}$ is the natural torus action 
of ${\mathbb T} := {\bold R}/T_0{\mathbb Z}$ induced by $\sigma$ {\rm (}$\sigma$ 
has the period $T_0$ by the assumption {\rm(ii))} with the natural quotient mapping 
$t \in {\bold R} \mapsto [t] \in {\mathbb T}$.  
\end{lemma} 

\medskip
This lemma enables us to compute the relative commutant 
$$
\left(\widetilde{{\mathcal N}(\pi)}\right)' \cap \widetilde{{\mathcal M}(\pi)}.
$$ 
(One should here note that this relative commutant involves continuous 
crossed-products so that the ``Fourier series expansion" cannot be used 
to deal with generic elements in any sense !) 

\medskip
We have 
$$
\begin{aligned}
&\left(M^{\Gamma}\rtimes_{\sigma}{\bold R}\right)' \cap 
\left({\bold C}1\rtimes{\bold R}\right) \\
&\cong 
\left(\left(M^{\Gamma}\rtimes_{\dot{\sigma}}{\mathbb T}\right)\otimes
L^{\infty}[0,-\log q^2)\right)' \cap 
\left(\left({\bold C}1\rtimes{\mathbb T}\right)\otimes
L^{\infty}[0,-\log q^2)\right) \\
&= {\bold C}1\otimes L^{\infty}[0,-\log q^2) 
\end{aligned}
$$
(the first isomorphism given by $\Phi$), and thus  
$$
\left(M^{\Gamma}\rtimes_{\sigma}{\bold R}\right)' \cap 
\left({\bold C}1\rtimes{\bold R}\right) 
= \left\{\lambda^{\sigma}(T_0) = 1\otimes\lambda_{T_0}\right\}''
$$
since 
$$
\left\{\Phi\left(\lambda^{\sigma}(T_0)\right)\right\}'' = 
\left\{1\otimes m_{e^{iT_0(\ \cdot\ )}}\right\}'' = 
{\bold C}1\otimes L^{\infty}[0,-\log q^2). 
$$
Therefore, we obtain 
$$
\begin{aligned}
&\left(\widetilde{{\mathcal N}(\pi)}\right)' \cap \widetilde{{\mathcal M}(\pi)} \\
&\subseteq 
\left({\bold C}1\otimes Q\otimes{\bold C}1\right)' \cap 
\left(B(V_{\pi})\otimes M)^{\text{\rm Ad}\pi\otimes\Gamma}\otimes B(L^2({\bold R})\right) \\
&= 
\left(B(V_{\pi})\otimes\left(Q'\cap M\right)\right)^{\text{\rm Ad}\pi\otimes\Gamma}\otimes
B(L^2({\bold R})) \quad \text{($Q$ sits in $M^{\Gamma}$)} \\
&= B(V_{\pi})^{\text{\rm Ad}\pi}\otimes{\bold C}1\otimes B(L^2({\bold R})) \quad \quad 
\text{(by the assumption (iii))} \\ 
&= \text{\rm End}(\pi)\otimes{\bold C}1\otimes B(L^2({\bold R})).
\end{aligned}
$$ 
Since 
$\sigma_t^{\psi\circ E_{\pi}} = 
\left(\text{Ad}F_{\pi}^{it}\otimes\sigma_t^{\varphi}\right)|_{{\mathcal M}(\pi)}$ and 
$F_{\pi} \in \text{\rm End}(\pi)'$ (see Eq.(3.9) in p.47, lines 3--6 from the top of p.49 of 
\cite{ueda:index4}), we have 
$$
\pi_{\sigma}(x\otimes1) = x\otimes1\otimes1 \quad \text{for all $x \in \text{\rm End}(\pi)$}, 
$$
and hence
$$
\left(\widetilde{{\mathcal N}(\pi)}\right)' \cap \widetilde{{\mathcal M}(\pi)} \subseteq 
\text{\rm End}(\pi)\otimes{\bold C}1\otimes\{\lambda_t : t \in {\bold R}\}''
$$
because $\{\lambda_t : t \in {\bold R}\}''$ is a MASA in $B(L^2({\bold R}))$. Here, $\lambda_t$ 
means the (left) regular representation of ${\bold R}$. Therefore, we get 
$$
\begin{aligned}
\left(\widetilde{{\mathcal N}(\pi)}\right)' \cap \widetilde{{\mathcal M}(\pi)} 
&= 
\left(\widetilde{{\mathcal N}(\pi)}\right)' \cap 
\left(\text{\rm End}(\pi)\otimes{\bold C}1\otimes\left\{\lambda_t : t \in {\bold R}\right\}''\right) \\
&= 
\left({\bold C}1\otimes\left(M^{\Gamma}\rtimes_{\sigma}{\bold R}\right)\right)' \cap 
\left(\text{\rm End}(\pi)\otimes\left({\bold C}\rtimes{\bold R}\right)\right) \\
&= 
\text{\rm End}(\pi)\otimes{\bold C}1\otimes\left\{\lambda_{T_0}\right\}'' \\
&\left(\cong \text{\rm End}(\pi)\otimes L^{\infty}[0,-\log q^2)\right). 
\end{aligned}
$$

\medskip
Set 
$$
{\mathcal Z} := 
{\mathcal Z}\left(\left(\widetilde{{\mathcal N}(\pi)} \right)' \cap \widetilde{{\mathcal M}(\pi)} \right) 
= {\mathcal Z}\left( {\rm End}(\pi) \right) \otimes{\bold C}1\otimes\left\{\lambda_{T_0}\right\}''.
$$
We have  
$$
\theta^{{\mathcal M}(\pi)}_t\left(x\otimes1\otimes\lambda_{T_0}\right) = 
x\otimes1\otimes\left(e^{-itT_0}\lambda_{T_0}\right), \quad 
x \in {\mathcal Z}\left(\text{\rm End}(\pi)\right),
$$
and hence, via the identification due to Lemma 1 (or $\Phi$),  
$$
\theta_t^{{\mathcal M}(\pi)} : m_{e^{iT_0(\ \cdot\ )}}\ \longmapsto\  
m_{e^{iT_0(\ \cdot\ -t)}} = m_{(e^{(iT_0(\ \cdot\ ))}\circ F_{-t})}, 
$$
where $F_t$ is the translation by $t$ (mod:$-\log q^2$). Thus, we get the isomorphism  
$$
\Psi : 
\left({\mathcal Z}, \theta^{{\mathcal M}(\pi)}_t|_{{\mathcal Z}}\right) \cong 
\left({\mathcal Z}\left(\text{\rm End}(\pi)\right)\otimes L^{\infty}[0,-\log q^2), 
\text{Id}\otimes F_t^*\right). 
$$

\medskip
What we want to do at first is to describe the relative position of the two subalgebras:  
$$
{\mathcal Z}\left(\widetilde{{\mathcal M}(\pi)}\right)\ \text{and}\  
{\mathcal Z}\left(\widetilde{{\mathcal N}(\pi)}\right) 
$$
inside ${\mathcal Z} = 
{\mathcal Z}\left(\left(\widetilde{{\mathcal N}(\pi)}\right)' \cap \widetilde{{\mathcal M}(\pi)}\right)$ by giving an explicit expression of the factor maps 
$$
X_{{\mathcal M}(\pi)} \overset{\pi^{{\mathcal M}(\pi)}}\longleftarrow 
X \overset{\pi^{{\mathcal N}(\pi)}}\longrightarrow X_{{\mathcal N}(\pi)} 
$$
with the point-realizations: 
$$
{\mathcal Z} = L^{\infty}(X), \quad 
{\mathcal Z}\left(\widetilde{{\mathcal M}(\pi)}\right) = L^{\infty}\left(X_{{\mathcal M}(\pi)}\right), \quad 
{\mathcal Z}\left(\widetilde{{\mathcal N}(\pi)}\right) = L^{\infty}\left(X_{{\mathcal N}(\pi)}\right).
$$

\medskip
The lemma below seems to be quite standard, but we will give the proof for the sake of completeness. 
 
\begin{lemma} We have, via the $\Psi$, the following description{\rm :}  
$$
\begin{aligned}
{\mathcal Z} \quad \ 
&= \
{\mathcal Z}\left(\text{\rm End}(\pi)\right)\otimes L^{\infty}[0,-\log q^2), \\
{\mathcal Z}\left(\widetilde{{\mathcal M}(\pi)}\right) \ 
&= \quad \left\{F_{\pi}^{-iT_0}\otimes m_{e^{iT_0(\ \cdot\ )}}\right\}'', \\
{\mathcal Z}\left(\widetilde{{\mathcal N}(\pi)}\right) \ 
&= \quad {\bold C}1\otimes L^{\infty}[0,-\log q^2). 
\end{aligned}
$$
\end{lemma}
\begin{proof} Only the non-trivial part is to compute 
${\mathcal Z}\left(\widetilde{{\mathcal M}(\pi)}\right)$. 
Since $\sigma_{T_0}^{\psi\circ E_{\pi}} = \text{\rm Ad}\left(F_{\pi}^{iT_0}\otimes1\right)$, 
the $F_{\pi}^{iT_0}\otimes1$ is in 
${\mathcal Z}\left({\mathcal M}(\pi)_{\psi\circ E_{\pi}}\right)$ so that we can find 
an invertible element $h \in {\mathcal Z}\left({\mathcal M}(\pi)_{\psi\circ E_{\pi}}\right)$ 
with $F_{\pi}^{iT_0}\otimes1 = h^{iT_0}$. (Note that $h$ and $F_{\pi}\otimes1$ are different 
and agree with only in the $iT_0$-power !) Perturbing $\psi\circ E_{\pi}$ by $h^{-1}$, i.e., 
$\psi\circ E_{\pi}(h^{-1}\ \cdot\ ) = \psi\circ E_{\pi}(h^{-1/2}\ \cdot\ h^{-1/2})$, we get 
the new state $\phi$ on ${\mathcal M}(\pi)$ with the property 
$\sigma_{T_0}^{\phi} = \text{\rm Id}$. The unitary operator $U^h$ on 
$L^2\left({\bold R}, L^2\left({\mathcal M}(\pi)\right)\right)$ defined by 
$$
\left(U^h \xi\right)(s) := h^{-is}\xi(s), \quad 
\xi \in L^2\left({\bold R}, L^2\left({\mathcal M}(\pi)\right)\right)
$$
gives us the isomorphism $\text{\rm Ad}U^h : {\mathcal M}(\pi)\rtimes_{\sigma^{\phi}}{\bold R} 
\cong {\mathcal M}(\pi)\rtimes_{\sigma^{\psi\circ E_{\pi}}}{\bold R}$, which satisfies that 
$$
\begin{aligned}
\text{\rm Ad}U^h\left(\pi_{\sigma^{\phi}}(m)\right) 
&= \pi_{\sigma^{\psi\circ E_{\pi}}}(m), \\
\text{\rm Ad}U^h\left(\lambda^{\sigma^{\phi}}(t)\right) 
&= 
\left(h^{-it}\otimes1\right)\lambda^{\sigma^{\psi\circ E_{\pi}}}(t)
= 
\pi_{\sigma^{\psi\circ E_{\pi}}}\left(h^{-it}\right)\lambda^{\sigma^{\psi\circ E_{\pi}}}(t). 
\end{aligned}
$$
Hence, we have 
$$
\begin{aligned}
{\mathcal Z}\left(\widetilde{{\mathcal M}(\pi)}\right) &= 
{\mathcal Z}\left({\mathcal M}(\pi)\rtimes_{\sigma^{\psi\circ E_{\pi}}}{\bold R}\right) \\
&= 
\text{\rm Ad}U^h\left({\mathcal Z}\left({\mathcal M}(\pi)\rtimes_{\sigma^{\phi}}{\bold R}\right)\right) \\
&= 
\text{\rm Ad}U^h\left(\left\{\lambda^{\sigma^{\phi}}(T_0)\right\}''\right) \\
&\phantom{aaaa} 
\text{(since $\sigma_{T_0}^{\phi} = \text{\rm Id}$, see around 
Lemma \ref{haagerupstormerlemma})}\\
&= 
\left\{\pi_{\sigma^{\psi\circ E_{\pi}}}\left(F_{\pi}^{-iT_0}\otimes1\right)
\lambda^{\sigma^{\psi\circ E_{\pi}}}(T_0)\right\}'' 
\quad \text{(since $F_{\pi}^{iT_0}\otimes1 = h^{iT_0}$)} \\
&= 
\left\{F_{\pi}^{-iT_0}\otimes1\otimes\lambda_{T_0}\right\}''. 
\end{aligned}
$$
Therefore, we get the desired identification via the isomorphism $\Psi$.   
\end{proof}

\medskip
It is known that the irreducible corepresentations are indexed by the spins, i.e., 
the half-integers so that we can associate the half-integer $\ell(\pi_j) \in 
\frac{1}{2}{\mathbb N}_0$ to any $\pi_j$, and we easily observe 
$$
F_{\pi}^{-iT_0} = F_{\pi}^{iT_0} = \sum (-1)^{2\ell(\pi_j)} p_j \in 
{\bold C}p_1\oplus\cdots\oplus{\bold C}p_n = {\mathcal Z}\left(\text{\rm End}(\pi)\right)
$$
based on the explicit computation of the matrix elements of the spin $\ell$ corepresentations, 
see \cite{mmnnu}. The standard identification 
$$
{\mathcal Z}\left(\text{\rm End}(\pi)\right)\otimes L^{\infty}[0,-\log q^2) \cong  
L^{\infty}\left(\{1,\dots,n\}\times[0,-\log q^2)\right)
$$
sends the generators 
$$
F_{\pi}^{-iT_0}\otimes m_{e^{iT_0(\ \cdot\ )}}, \quad  1\otimes m_{e^{iT_0(\ \cdot\ )}}
$$
of ${\mathcal Z}\left(\widetilde{{\mathcal M}(\pi)}\right)$, 
${\mathcal Z}\left(\widetilde{{\mathcal N}(\pi)}\right)$ to the functions  
$$
\begin{aligned}
u^{{\mathcal M}(\pi)}(j,s) &:= (-1)^{2\ell(\pi_j)} e^{iT_0 s} = 
e^{iT_0(s - \ell(\pi_j) \log q^2)}, \\
u^{{\mathcal N}(\pi)}(j,s) &:= e^{iT_0 s}, 
\end{aligned}
$$
respectively. The mapping 
$$
\pi_0 : (j,s) \mapsto (j,s - \ell(\pi_j) \log q^2)
$$ defined on $\{1,\dots,n\}\times[0,-\log q^2)$ satisfies 
$$
u^{{\mathcal M}(\pi)}\circ\pi_0^{-1} = u^{{\mathcal N}(\pi)} \quad  
\text{or} \quad 
u^{{\mathcal M}(\pi)} = u^{{\mathcal N}(\pi)}\circ\pi_0 
$$ 
so that the induced automorphism $\Phi_0$ of ${\mathcal Z}$ of $\pi_0$ satisfies  
$$
\Phi_0\left({\mathcal Z}\left(\widetilde{{\mathcal M}(\pi)}\right)\right) =  
{\mathcal Z}\left(\widetilde{{\mathcal N}(\pi)}\right), \quad 
\left(\theta_t\big|_{\mathcal Z}\right)\circ\Phi_0 = 
\Phi_0\circ\left(\theta_t\big|_{\mathcal Z}\right).
$$ 
Since both ${\mathcal M}(\pi)$ and ${\mathcal N}(\pi)$ are factors of type III$_{q^2}$, 
the $X_{{\mathcal M}(\pi)}$, $X_{{\mathcal N}(\pi)}$ are the same 
space; 
$$
X_{{\mathcal M}(\pi)} = X_{{\mathcal N}(\pi)} = [0,-\log q^2),  
$$
and the identification here enables us to describe the factor map  
$$
\pi^{{\mathcal N}(\pi)} : 
X : = \{1,\dots,n\}\times[0,-\log q^2) \rightarrow X_{{\mathcal N}(\pi)} = [0,-\log q^2)
$$ 
by the standard one $(j,s) \mapsto s$ while the factor map 
$$
\pi^{{\mathcal M}(\pi)} : 
X := \{1,\dots,n\}\times[0,-\log q^2) \rightarrow X_{{\mathcal M}(\pi)} = [0,-\log q^2)
$$
by the deformation of the standard one by $\pi_0$; $\pi^{{\mathcal M}(\pi)} = \pi^{{\mathcal N}(\pi)}
\circ\pi_0$, that is, $(j,s) \mapsto s - \ell(\pi_j) \log q^2$ (mod: $-\log q^2$). Hence, conclude

\begin{theorem} 
The factor map structure attached to ${\mathcal M}(\pi) \overset{E_{\pi}}\supseteq {\mathcal N}(\pi)$ 
is  
$$
[0,-\log q^2) 
\overset{\pi^{{\mathcal M}(\pi)}}\longleftarrow 
\{1,\dots,n\}\times[0,-\log q^2) 
\overset{\pi^{{\mathcal N}(\pi)}}\longrightarrow 
[0,-\log q^2)
$$
with 
$$
\pi^{{\mathcal M}(\pi)}(j,s) = s -\ell(\pi_j) \log q^2\ \text{\rm (mod: $-\log q^2$)}, \quad 
\pi^{{\mathcal N}(\pi)}(j,s) = s. 
$$
Therefore, the inclusion is of essentially type {\rm II} if and only if the set 
$$
\left\{ \ell : \text{the spin $\ell$ corepresentation occurs in $\pi$} \right\} 
$$
consists entirely of integers or consists entirely of half-integers. 
\end{theorem}

\medskip
The next question from the view-point of Hamachi-Kosaki's analysis is to describe the 
following inclusion of von Neumann algebras (not necessary being factors) inside 
$\displaystyle{{\mathcal M}(\pi) \supseteq {\mathcal N}(\pi)}$: Letting 
$$
\widetilde{{\mathfrak A}} := \widetilde{{\mathcal M}(\pi)}\cap{\mathcal Z}' \supseteq 
\widetilde{{\mathfrak B}} := \widetilde{{\mathcal N}(\pi)}\vee{\mathcal Z}, 
$$
we set
$$
{\mathfrak A} := \widetilde{\mathfrak A}\rtimes_{\theta}{\bold R} \supseteq 
{\mathfrak B} := \widetilde{\mathfrak B}\rtimes_{\theta}{\bold R}
$$
with $\theta := \theta^{{\mathcal M}(\pi)}$. We first easily observe that  
$$
\begin{aligned}
{\mathcal Z}\left(\widetilde{{\mathfrak A}}\right) &= 
{\mathcal Z}\left(\widetilde{{\mathfrak B}}\right) = 
{\mathcal Z}, \\
{\mathcal Z}\left({\mathfrak A}\right) &= 
{\mathcal Z}\left({\mathfrak B}\right) = 
{\mathcal Z}^{\theta} = {\mathcal Z}\left(\text{\rm End}(\pi)\right)\otimes{\bold C}1\otimes{\bold C1} 
\end{aligned}
$$
since $F_t$ is a transitive flow on $[0,-\log q^2)$. What we will do is to give a complete 
description of this inclusion ${\mathfrak A} \supseteq {\mathfrak B}$ in terms of the irreducible decomposition of $\pi$. 

Consider the conditional expectation 
$$
\widehat{G} : 
\widetilde{{\mathcal M}(\pi)} \rightarrow  \widetilde{\mathfrak A} = 
\widetilde{{\mathcal M}(\pi)} \cap 
\left(\pi_{\sigma}\left({\mathcal Z}\left(\text{\rm End}(\pi)\right)
\otimes{\bold C}1\right)\right)'
$$
given by 
$$
\widehat{G}(X) := 
\int_{{\mathcal U}\left({\mathcal Z}\left(\text{\rm End}(\pi)\right)\right)} 
\pi_{\sigma}(u\otimes1) X \pi_{\sigma}(u\otimes1)^* du, \quad X \in \widetilde{{\mathcal M}(\pi)}, 
$$
with the modular automorphism $\sigma_t := \sigma^{\psi\circ E_{\pi}}_t$. 
For any finite sum 
$$
X = \sum_j\pi_{\sigma}(m(t_j)) \lambda^{\sigma}(t_j) \in 
\widetilde{{\mathcal M}(\pi)}
$$ 
we have 
$$
\begin{aligned}
\widehat{G}(X) &= 
\sum_j \widehat{G}\left(\pi_{\sigma}\left(m(t_j)\right)\lambda^{\sigma}(t_j)\right) \\
&= 
\sum_j 
\int_{{\mathcal U}\left({\mathcal Z}\left(\text{\rm End}(\pi)\right)\right)} 
\pi_{\sigma}\left((u\otimes1)m(t_j)(u\otimes1)^*\right)\lambda^{\sigma}(t_j) du \\
&= 
\sum_j 
\pi_{\sigma}\left(\int_{{\mathcal U}\left({\mathcal Z}\left(\text{\rm End}(\pi)\right)\right)} 
(u\otimes1)m(t_j)(u\otimes1)^*du\right)\lambda^{\sigma}(t_j).
\end{aligned}
$$
Here the second equality comes from the fact that the modular automorphism 
$\sigma_t$ acts on ${\mathcal Z}(\text{\rm End}(\pi))\otimes{\bold C}1$ trivially. 
Note that all the 
$$
\int_{{\mathcal U}\left({\mathcal Z}\left(\text{\rm End}(\pi)\right)\right)} 
(u\otimes1)m(t_j)(u\otimes1)^*du 
$$
are in the relative commutant 
$$
\left(B(V_{\pi})\otimes M\right)^{\text{Ad}\pi\otimes\Gamma} \cap 
\left({\mathcal Z}\left(\text{\rm End}(\pi)\right)\otimes{\bold C}1\right)' 
= \left({\mathcal Z}\left(\text{\rm End}(\pi)\right)'\otimes M\right)^{\text{Ad}\pi\otimes\Gamma}. 
$$
Here we used that ${\mathcal Z}\left(\text{\rm End}(\pi)\right)\otimes{\bold C}1$ sits in 
$\left(B(V_{\pi})\otimes M\right)^{\text{Ad}\pi\otimes\Gamma}$. 
Thus, by the appropriate continuity of $\widehat{G}$, 
$$
\begin{aligned}
\widetilde{\mathfrak A} 
&= 
\widetilde{{\mathcal M}(\pi)} \cap 
\left(\pi_{\sigma}\left({\mathcal Z}\left(\text{\rm End}(\pi)\right)
\otimes{\bold C}1\right)\right)' \\
&= 
\widehat{G}\left(\widetilde{{\mathcal M}(\pi)}\right) \subseteq 
\left({\mathcal Z}\left(\text{\rm End}(\pi)\right)'\otimes M\right)^{\text{\rm Ad}\pi\otimes\Gamma}
\rtimes_{\sigma}{\bold R}, 
\end{aligned}
$$
and hence 
$$
\widetilde{\mathfrak A} = 
\left({\mathcal Z}\left(\text{\rm End}(\pi)\right)'\otimes M\right)^{\text{\rm Ad}\pi\otimes\Gamma}
\rtimes_{\sigma}{\bold R}. 
$$
It is clear that 
$$
\begin{aligned}
\widetilde{{\mathfrak B}} &= \widetilde{{\mathcal N}(\pi)} \vee {\mathcal Z} \\
&= 
\left( {\bold C}1 \otimes \left( M^{\Gamma} \rtimes_{\sigma} {\bold R} \right) \right) \vee 
\left({\mathcal Z}\left(\text{\rm End}(\pi)\right)
\otimes{\bold C}1\otimes\left\{\lambda_{T_0}\right\}'' \right) \\
&= 
{\mathcal Z}\left(\text{\rm End}(\pi)\right)\otimes\left(M^{\Gamma}\rtimes_{\sigma}{\bold R}\right). 
\end{aligned}
$$
Via the Takesaki duality, we get 
$$
{\mathfrak A} = 
\left({\mathcal Z}\left(\text{\rm End}(\pi)\right)'\otimes M\right)^{\text{\rm Ad}\pi\otimes\Gamma}, 
\quad 
{\mathfrak B} = 
{\mathcal Z}\left(\text{\rm End}(\pi)\right)\otimes M^{\Gamma}. 
$$  

We will next give the description of the decomposition of $E_{\pi} : {\mathcal M}(\pi) \rightarrow 
{\mathcal N}(\pi)$ into 
$$
{\mathcal M}(\pi) \overset{G}\longrightarrow {\mathfrak A} \overset{H}\longrightarrow 
{\mathfrak B} \overset{F}\longrightarrow {\mathcal N}(\pi)
$$
with $E_{\pi} = F\circ H\circ G$. Define the conditional expectations 
$$
\begin{aligned}
F &: 
{\mathfrak B} = {\mathcal Z}\left(\text{\rm End}(\pi)\right)\otimes M^{\Gamma} \rightarrow 
{\mathcal N}(\pi) = {\bold C}1\otimes M^{\Gamma}; \\
H &: 
{\mathfrak A} = 
\left({\mathcal Z}\left(\text{\rm End}(\pi)\right)'\otimes M\right)^{\text{\rm Ad}\pi\otimes\Gamma} 
\rightarrow 
{\mathfrak B} = {\mathcal Z}\left(\text{\rm End}(\pi)\right)\otimes M^{\Gamma}; \\
G &: 
{\mathcal M}(\pi) = 
\left(B(V_{\pi})\otimes M\right)^{\text{\rm Ad}\pi\otimes\Gamma} \rightarrow 
{\mathfrak A} = 
\left({\mathcal Z}\left(\text{\rm End}(\pi)\right)'\otimes M\right)^{\text{\rm Ad}\pi\otimes\Gamma} 
\end{aligned}
$$
by 
$$
\begin{aligned}
F(X) &:= \left(\tau_q^{(\pi,-)}\otimes\text{\rm Id}\right)(X), \quad X \in {\mathfrak B}; \\
H(X) &:= 
\sum_{j=1}^n \frac{1}{\tau_q^{(\pi,-)}(p_j)} \left(\tau_q^{(\pi,-)}\otimes\text{\rm Id}\right)
\left(X(p_j\otimes1)\right)p_j\otimes1, \quad X \in {\mathfrak A}; \\
G(X) &:= 
\int_{{\mathcal U}\left({\mathcal Z}\left(\text{\rm End}(\pi)\right)\right)} 
(u\otimes1) X (u^*\otimes1) du \\
&= \sum_{j=1}^n{}^{\oplus} (p_j\otimes1)X(p_j\otimes1), \quad X \in {\mathfrak B},  
\end{aligned}
$$   
where the $p_j$'s form the partition of unity consisting of minimal projections in 
${\mathcal Z}\left(\text{\rm End}(\pi)\right)$. Here, the last expression of $G$ is 
due to the simple fact that 
$$
\int_{\mathbb T} \int_{\mathbb T} z_1\bar{z}_2 d\mu_{\mathbb T}(z_1)d\mu_{\mathbb T}(z_2)  
= 0
$$ 
with the Haar probability measure $\mu_{\mathbb T}$. 
Since $\sigma_t = \sigma_t^{\psi\circ E_{\pi}}$ acts on 
${\mathcal Z}\left(\text{\rm End}(\pi)\right)\otimes{\bold C}1$ trivially as remarked before, 
we have, by a direct computation,
$$
\psi\circ F\circ H\circ G(X) = \psi\circ E_{\pi}(X), \quad X \in {\mathcal M}(\pi), 
$$
and hence $E_{\pi} = F\circ H\circ G$. One should note that $G$ agrees with $\widehat{G}$ 
introduced previously. 

\medskip
Therefore, putting the detailed data of  
$$
\pi = m_1\pi_1\oplus\cdots\oplus m_n\pi_n
$$
into the above computations we finally get the following descriptions: 
$$
\begin{aligned}
&\left({\mathfrak A} \supseteq {\mathfrak B}\right) \\
&= 
\sum_{j=1}^n{}^{\oplus} \left(
M_{m_j}({\bold C})\otimes\left(B(V_{\pi_j})\otimes M\right)^{\text{\rm Ad}\pi_j\otimes\Gamma} 
\supseteq  
{\bold C}1_{m_j}\otimes{\bold C}1_{V_{\pi_j}}\otimes M^{\Gamma}
\right) \\
&= 
\sum_{j=1}^n{}^{\oplus} \left( M_{m_j}({\bold C})\otimes{\mathcal M}(\pi_j) \supseteq 
{\bold C}1_{m_j}\otimes{\mathcal N}(\pi_j)
\right),  
\end{aligned}
$$
where the $\left({\mathcal M}(\pi_j) \supseteq {\mathcal N}(\pi_j)\right)$'s mean 
the Wassermann-type inclusions associated with the $\pi_j$'s. The conditional expectation 
$H : {\mathfrak A} \rightarrow {\mathfrak B}$ can be decomposed into 
$$
H_j := H|_{(p_j\otimes1)} : {\mathfrak A}_{(p_j\otimes1)} \rightarrow {\mathfrak B}_{(p_j\otimes1)}, 
$$ 
and we have 
$$
\begin{aligned}
H_j(\ \cdot\ ) &
= 
\frac{1}{\tau_q^{(\pi,-)}(p_j)} \left(\tau_q^{(\pi,-)}\otimes\text{\rm Id}\right)
\left(\ \cdot\ (p_j\otimes1)\right)(p_j\otimes1)  \\
&= 
\frac{\text{dim}_q(\pi)}{\text{\rm Tr}(1_{m_j}\otimes F_{\pi_j}^{-1})} \cdot 
\frac{1}{\text{dim}_q(\pi)} \left(\text{\rm Tr}
\left(1\otimes F_{\pi_j}^{-1}\ \cdot\ \right)\otimes\text{\rm Id}\right)(p_j\otimes1) \\
&= 
\frac{1}{m_j\cdot\text{\rm dim}_q(\pi_j)} \left(\text{\rm Tr}_{m_j}\otimes\text{Tr}_{V_{\pi_j}}
\left(F_{\pi_j}^{-1}\ \cdot\ \right)\otimes\text{\rm Id}\right)(p_j\otimes1) \\
&= 
\left(\left(\frac{1}{m_j}\text{\rm Tr}_{m_j}\right)\otimes
\left(\frac{1}{\text{\rm dim}_q(\pi_j)}\text{Tr}_{V_{\pi_j}}
\left(F_{\pi_j}^{-1}\ \cdot\ \right)\right)\otimes
\text{\rm Id}\right)(p_j\otimes1) \\
&= 
\left(\tau_{m_j}\otimes\left(\tau_q^{(\pi_j,-)}\otimes\text{\rm Id}\right)\right)|_
{M_{m_j}({\bold C})\otimes\left(B(V_{\pi_j})\otimes M\right)^{\text{\rm Ad}\pi_j\otimes\Gamma}},  
\end{aligned}
$$
where we denote the non-normalized usual trace by $\text{\rm Tr}$ and the normalized one by $\tau$. 
Thus the conditional expectation $H : {\mathfrak A} \rightarrow {\mathfrak B}$ has, via the above 
identification, the following central decomposition: 
$$
H = \sum_{j=1}^n{}^{\oplus}\ \tau_{m_j}\otimes E_{\pi_j},    
$$
where the $E_{\pi_j}$'s are the conditional expectations of the Wassermann-type inclusions  
$\left({\mathcal M}(\pi_j) \supseteq {\mathcal N}(\pi_j)\right)$'s. Summing up the discussions, 
we conclude   

\medskip
\begin{theorem} The three-step inclusions 
$$
{\mathcal M}(\pi)\ \overset{G}\supseteq\ {\mathfrak A}\ \overset{H}\supseteq\ {\mathfrak B}\ 
\overset{F}\supseteq\ {\mathcal N}(\pi)
$$
have the following description: 
$$
\left({\mathfrak A} \overset{H}\supseteq {\mathfrak B}\right) \cong  
\sum_{j=1}^n{}^{\oplus}  
\left( M_{m_j}({\bold C})\otimes{\mathcal M}(\pi_j) \overset{\tau_{m_j}\otimes E_{\pi_j}}\supseteq 
{\bold C}1_{m_j}\otimes{\mathcal N}(\pi_j) \right)  
$$
and the conditional expectations $G$, $F$ are given by 
$$
G(\ \cdot\ ) = \sum_{j=1}^n{}^{\oplus} (p_j\otimes1)\left(\ \cdot\ \right)(p_j\otimes1), \quad
F = E_{\pi}\big|_{\mathfrak B},  
$$
where the $p_j$'s mean the partition of unity consisting of minimal projections in 
${\mathcal Z}\left(\text{\rm End}(\pi)\right)$.  
\end{theorem} 

\medskip
One should note that the factors ${\mathcal M}(\pi_j)$, ${\mathcal N}(\pi_j)$, etc., appearing in 
the above expression are all isomorphic to $M^{\Gamma}$ by the construction itself. 

\section{Comments and Questions.} 
\label{comments}

(1) {\it The existence question of AFD-minimal actions}: We provide, in 
\cite{ueda:minimalaction,ueda:index4}, examples of minimal actions 
of compact non-Kac quantum groups. However, their constructions involve  
the free product construction so that the factors on which the quantum groups in question act 
do  never be AFD. Hence, the most important question in the subject matter is: 
Does a compact non-Kac quantum group possesses minimal actions on 
AFD-factors ? If it was affirmative, we would like to ask further: What kind of compact 
non-Kac quantum groups does allow minimal actions on AFD factors ?  
 
\medskip\noindent
(2) {\it Bernoulli shift actions of discrete quantum groups}: In a more or less connection to 
the above-mentioned questions, we would like to pose the research project to find true analog(s?) 
of (non-commutative) Bernoulli shift actions of discrete quantum groups (or equivalently 
of the duals of compact (quantum) groups, e.g. $\widehat{{\rm SU}_q(2)}$).  The duality given in  
Theorem \ref{hiai-ueda-duality} might be helpful to try this research project. Indeed, as was 
commented after Theorem \ref{hiai-ueda-duality} the theorem explains philosophically that 
our examples of minimal actions of ${\rm SU}_q(N)$ are thought of as the dual actions of 
``free Bernoulli shift actions" of the dual $\widehat{{\rm SU}_q(N)}$. Thus, the first attempt 
seems to try to find an explicit construction of these ``free Bernoulli shifts" out of the dual 
quantum group $\widehat{{\rm SU}_q(N)}$ directly.  However, the most interesting question 
on this topic is what should be the discrete quantum group version of non-commutative 
Bernoulli $G$-shifts (with discrete group $G$) based on the infinite tensor product construction, 
which might help to solve the questions in (1).  
    
\medskip\noindent
(3) {\it Around free Bernoulli schemes}: Let $p = (p_1,\dots,p_n)$ be a probability vector, 
and there are two kinds of free analogs of the classical Bernoulli scheme associated with $p$ -- 
the ``free commutative Bernoulli scheme" ${\rm FCBS}(p) = \left(L({\mathbb F}_{\infty}), 
\gamma_p\right)$ and ``free non-commutative Bernoulli scheme" ${\rm FNBS}(p) = 
\left(L({\mathbb F}_{\infty}), \sigma_p\right)$, where $\gamma_p$ and $\sigma_p$ are both 
aperiodic, ergodic automorphisms on $L({\mathbb F}_{\infty})$. See \cite{hiai-ueda} 
for their constructions. These free Bernoulli schemes need new dynamical invariants like entropies. 
Indeed, any kind of non-commutative dynamical entropy cannot be used to distinguish them. In fact, 
it is known that they always have Connes-St{\o}rmer's entropy (see St{\o}rmer\cite{stormer}) 
$H^{\rm CS}_{\tau_{{\mathbb F}_{\infty}}}(\gamma_p) = 
H^{\rm CS}_{\tau_{{\mathbb F}_{\infty}}}(\sigma_p) = 0$ (for any $p$),  
and also we can easily show that the perturbation theoretic dynamical entropy (see 
\cite{voiculescu:perturbationtheoreticentropy}) 
$H_P(\gamma_p) = H_P(\sigma_p) = \infty$ (for any $p$). 
However, we have very small facts on 
their crossed-products (see Proposition 14, 15 of \cite{hiai-ueda}), which in particular 
says that, if all the interpolated free group factors were mutually non-isomorphic, then we would 
have some cocycle conjugacy classification results for free Bernoulli schemes. 
Moreover, K.~Dykema \cite{dykema:comment} pointed out us that the questions of cocycle 
conjugacy and of conjugacy for ${\mathbb Z}$-free shifts considered in Proposition 12 of 
\cite{hiai-ueda} are both equivalent to the isomorphism question of their crossed-products, 
which can be easily derived from Theorem \ref{hiai-ueda-duality}. Thus, from the view-point of 
the famous isomorphism question of free group factors, it might be very interesting to seek for 
a new kind of non-commutative dynamical entropy fit for a suitable class of 
II$_1$-factors consisting of all free group factors, and the dynamical free entropy dimension 
(see \S\S7.2 of \cite{voiculescu:freeentropy2}) 
seems to work for those ${\mathbb Z}$-free shifts. (Unfortunately, its computation is probably  
quite difficult because of the same trouble (the so-called ``semi-continuity problem") for 
free entropy !)

\medskip\noindent
(4) {\it Some brief comments to} \S\ref{wassermannconstruction} {\it and} \S\ref{typeIIIsubfactoranalysis}: 
(4-1) It is possible to consider other quantum groups instead of ${\rm SU}_q(2)$. Here, we will briefly 
explain what phenomenon occurs in the ${\rm SO}_q(3)$-case. The quantum ${\rm SO}_q(3)$ 
($0 < q < 1$) is defined as the subalgebra of $L^{\infty}({\rm SU}_q(2))$ generated by the unitary 
$u(\pi_1)$ of the spin $1$ (irreducible, unitary) corepresentation $\pi_1$ of ${\rm SU}_q(2)$ with the 
same Hopf algebra structure, and it is known that the set of all nonequivalent irreducible corepresentations 
can be chosen as $\{\pi_{\ell}\}_{\ell \in {\mathbb N}_0}$. See \cite{podles:SO_q} for details. 
In this case, we can still prove that the fixed-point algebra $M^{\Gamma}$ (with minimal action $\Gamma$ 
defined in the same manner as in Theorem \ref{dima-ueda-maintheorem}) is of type III$_{q^2}$ (one of its 
proofs uses the explicit description of the fixed-point algebra by creation operators given in the Appendix I 
of \cite{dima-ueda}), and thus all the conditions given in the beginning of \S\ref{typeIIIsubfactoranalysis} 
are still satisfied, and all the results there are of course valid in this case. The most important thing in this 
case is the fact that the resulting inclusions always become of essentially type II because of all the irreducible 
corepresentations are indexed by only the integers.      

(4-2) Soon after the appearance of our work \cite{dima-ueda}, F.~R{\u{a}}dulescu \cite{radulescu:index4} 
gave another proof (or construction) to Theorem \ref{dima-ueda-maintheorem}. His subfactors are constructed 
by a similar (but not the exactly same) way as in \cite{popa:index4}. However, he also used a (minimal) 
action of ${\rm SU}_q(2)$ (or we should say a minimal action of ${\rm SO}_q(2)$) to identify with 
$L({\mathbb F}_{\infty})$. What we would like here to point out is the rather simple fact that his subfactors 
are also Wassermann-type. In fact,  one can easily see that the subfactors are of the form 
${\mathcal N} \supseteq {\mathcal N}_{-1}$, where ${\mathcal N}_{-1}$ is the downward basic extension 
of the Wassermann-type inclusion ${\mathcal M} \supseteq {\mathcal N}$ constructed from 
the free product action $\Gamma_{1/2} : M \rightarrow M\otimes A$: 
$$
\Gamma_{1/2} := (\text{Id}_Q \otimes 1_A) * \text{\rm Ad}\pi_{1/2}
$$
with $M := Q * B(V_{\pi_{1/2}})$.

\end{document}